\documentclass[a4paper,12pt]{amsart}
\usepackage{amssymb}
\usepackage{ifthen}
\usepackage{graphicx}
\usepackage{float}
\usepackage{caption}
\usepackage{subcaption}
\usepackage{cite}
\usepackage{amsfonts}
\usepackage{amscd}
\usepackage{amsxtra}
\usepackage{color}

\newtheorem{thm}{Theorem}
\newtheorem{cor}{Corollary}
\newtheorem{lem}{Lemma}

\newtheorem{rem}{Remark}
\newtheorem{example}{Example}
\newtheorem{defn}{Definition}
\newtheorem{prob}{Problem}

\newtheorem{conj}{Conjecture}
\theoremstyle{definition}

\newcounter {own}
\def\theown {\thesection  .\arabic{own}}

\newenvironment{pf}[1][]{%
 \vskip 3mm
 \noindent
 \ifthenelse{\equal{#1}{}}%
  {{\slshape Proof. }}%
  {{\slshape #1.} }%
 }%
{\qed\bigskip}

\newcounter{alphabet}
\newcounter{tmp}
\newenvironment{Thm}[1][]{\refstepcounter{alphabet}%
\bigskip%
\noindent%
{\bf Theorem \Alph{alphabet}}%
\ifthenelse{\equal{#1}{}}{}{ (#1)}%
{\bf .} \itshape}{\vskip 8pt}

\makeatletter
\newcommand{\Ref}[1]{\@ifundefined{r@#1}{}{\setcounter{tmp}{\ref{#1}}\Alph{tmp}}}
\makeatother

\newenvironment{Lem}[1][]{\refstepcounter{alphabet}%
\bigskip%
\noindent%
{\bf Lemma \Alph{alphabet}}%
{\bf .} \itshape}{\vskip 8pt}

\newcommand{\IR}{{\mathbb R}}
\newcommand{\IN}{{\mathbb N}}

\newcommand{\ID}{{\mathbb D}}

\def\be{\begin{equation}}
\def\ee{\end{equation}}

\newcommand{\bee}{\begin{enumerate}}
\newcommand{\eee}{\end{enumerate}}

\newcommand{\blem}{\begin{lem}}
\newcommand{\elem}{\end{lem}}
\newcommand{\bthm}{\begin{thm}}
\newcommand{\ethm}{\end{thm}}
\newcommand{\bcor}{\begin{cor}}
\newcommand{\ecor}{\end{cor}}
\newcommand{\beg}{\begin{example}}
\newcommand{\eeg}{\end{example}}
\newcommand{\begs}{\begin{examples}}
\newcommand{\eegs}{\end{examples}}
\newcommand{\bdefe}{\begin{defn}}
\newcommand{\edefe}{\end{defn}}
\newcommand{\bprob}{\begin{prob}}
\newcommand{\eprob}{\end{prob}}
\newcommand{\bei}{\begin{itemize}}
\newcommand{\eei}{\end{itemize}}

\newcommand{\bcon}{\begin{conj}}
\newcommand{\econ}{\end{conj}}
\newcommand{\bcons}{\begin{conjs}}
\newcommand{\econs}{\end{conjs}}
\newcommand{\bprop}{\begin{propo}}
\newcommand{\eprop}{\end{propo}}
\newcommand{\br}{\begin{rem}}
\newcommand{\er}{\end{rem}}
\newcommand{\brs}{\begin{rems}}
\newcommand{\ers}{\end{rems}}
\newcommand{\bo}{\begin{obser}}
\newcommand{\eo}{\end{obser}}
\newcommand{\bos}{\begin{obsers}}
\newcommand{\eos}{\end{obsers}}
\newcommand{\bpf}{\begin{pf}}
\newcommand{\epf}{\end{pf}}
\newcommand{\ba}{\begin{array}}
\newcommand{\ea}{\end{array}}
\newcommand{\beq}{\begin{eqnarray}}
\newcommand{\beqq}{\begin{eqnarray*}}
\newcommand{\eeq}{\end{eqnarray}}
\newcommand{\eeqq}{\end{eqnarray*}}

\newcommand{\ds}{\displaystyle}

\newcounter{minutes}
\divide\time by 60
\newcounter{hours}
\multiply\time by 60 \addtocounter{minutes}{-\time}

\begin{document}
\bibliographystyle{amsplain}
\title[Absolutely Convex,  Uniformly Starlike and Uniformly Convex Harmonic Mappings]
{Absolutely Convex,  Uniformly Starlike and Uniformly Convex Harmonic Mappings}

\thanks{
File:~\jobname .tex,
          printed: \number\year-\number\month-\number\day,
          \thehours.\ifnum\theminutes<10{0}\fi\theminutes}

\author{S. Ponnusamy
}
\address{S. Ponnusamy, and A. Sairam Kaliraj,
Indian Statistical Institute (ISI), Chennai Centre, SETS, MGR Knowledge City, CIT
Campus, Taramani, Chennai 600 113, India. }
\email{samy@isichennai.res.in, sairamkaliraj@gmail.com}

\author{A. Sairam Kaliraj}

\author{Victor V. Starkov}
\address{V.V. Starkov, Petrozavodsk State University, 33, Lenin Str., 185910, Petrozavodsk, Republic of Karelia, Russia.}
\email{vstarv@list.ru}

\keywords{Harmonic univalent, uniformly starlike, uniformly convex, absolutely convex,
growth theorem, covering theorem, coefficient bounds, Bernardi type operator, Gaussian hypergeometric functions.\\
$
^\dagger$ {\tt
The first author is currently on leave from the IIT Madras, India.
}
}

\date{\today  
~File:~pss3.tex}

\begin{abstract}
In this paper,  we prove necessary and sufficient conditions for a sense-preserving harmonic function to be absolutely convex in the open
unit disk. We also estimate the coefficient bound and obtain growth, covering and area theorems for absolutely convex harmonic mappings.
A natural generalization of the classical Bernardi type operator for harmonic functions is considered and its connection
between certain classes of uniformly starlike harmonic functions and uniformly convex harmonic functions is also investigated.
At the end, as applications, we present a number of results connected with hypergeometric and polylogarithm functions.
\end{abstract}
\thanks{ }

\maketitle
\pagestyle{myheadings}
\markboth{S.Ponnusamy, A. Sairam Kaliraj and V.V. Starkov}{Absolutely Convex and Uniformly Convex Harmonic Mappings}

\section{Introduction}
\setcounter{equation}{0}

Let $\mathcal{H}$ denote the class of all complex-valued harmonic functions $f=h+\overline{g}$ in $\mathbb{D} = \{z:\, |z|<1\}$, where $h$ and $g$ are analytic in $\mathbb{D}$ and normalized such that
\begin{equation}\label{intro1}
f(z)= h(z)+\overline{g(z)} = z + \sum_{n=2}^{\infty} a_nz^n  +\overline{\sum_{n=1}^{\infty} b_nz^n}.
 \end{equation}
The Jacobian of $f$ is given by $J_f(z) = \left|f_z(z)\right|^2 - \left|f_{\overline{z}}(z)\right|^2$ and we say that $f$ is sense-preserving in $\ID$ if
$J_f(z)>0$ in $\ID$.
Let $\mathcal{S}_{H}$ denote the subclass of $\mathcal{H}$ consisting of sense-preserving univalent functions in $\mathbb{D}$. Also, let $\mathcal{S}^0_{H} = \{ f \in \mathcal{S}_{H}: b_1=\overline{f_{\overline{z}}(0)}=0\}$ and $\mathcal{S}=\{ f = h + \overline{g}\in \mathcal{S}_{H}:\, g\equiv0\}$, the classical family of univalent analytic functions on $\ID$. The family $\mathcal{S}^0_{H}$ is known to be compact and normal, whereas $\mathcal{S}_{H}$ is normal but not compact. For discussion on this and several other geometric subclasses, see \cite{Clunie-Small-84, Duren:Harmonic,SaRa2013}. One of the important geometric subclasses is the class $\mathcal{K}_{H}$ of sense-preserving univalent harmonic mappings $f$ of the form \eqref{intro1} such that $f(\ID)$ is a convex domain. We may now set
$$\mathcal{K} =\{f= h + \overline{g}\in\mathcal{K}_H:\, g(z)\equiv 0 \}.
$$
For $f \in \mathcal{K}$, $f(\ID_r)$ is also convex for each $r \in (0, 1)$, where $\ID_r:=\{z:\, |z|<r\}$. In other words, we say that each $f \in \mathcal{K}$ is
fully convex in $\ID$. Unfortunately, this property does not carry over to functions in
$\mathcal{K}_H$. Consequently, it is natural to introduce and study a subclass of $\mathcal{K}_H$ having analogous property.
In \cite{Duren_Osgood} Duren et al. introduced the subclass $\mathcal{FK}_{H}$ of fully convex harmonic mappings. That is, $f \in \mathcal{FK}_{H}$ if and only if $f(\ID_r)$ is a convex domain for each $r \in (0, 1)$, where $\ID_r:=\{z:\, |z|<r\}$. Rad\'{o}-Kneser-Choquet theorem \cite[Section 3.1]{Duren:Harmonic} ensures that a fully convex harmonic mapping is necessarily univalent in $\mathbb{D}$.
It is worth recalling  that, if
$f \in \mathcal{K}$, then $f(\ID(a, r))$ is convex for all $a \in \ID$ and $r \in (0, 1)$ such that $|a|+r < 1$, where
$$\ID(a, r)=\{z\in\ID:\, |z-a| < r\}.
$$
We emphasize that this property does not hold for functions in $\mathcal{FK}_H$. This leads to the introduction of the following.

\bdefe A locally univalent harmonic function $f=h+\overline{g} \in \mathcal{H}$ is said to be absolutely convex in $\mathbb{D}$, symbolically denoted by $f \in \mathcal{AK}_H$, if $f(\ID(a, r))$ is convex for all $a \in \ID$ and for all $r \in (0, 1-|a|)$.
\edefe

We denote by $\mathcal{AK}^0_H=\{f = h + \overline{g}\in \mathcal{AK}_H:\, g'(0)=0 \}$ and
$$\widetilde{\mathcal{AK}_H}=\{ f \in \mathcal{AK}_H:\, f(\ID(a, r)) ~\mbox{is convex for all}~a \in \ID ~\mbox{and}~r \in (0, 1-|a|] \,\}.
$$


Denote by $\mathcal{UK}_{H}$ (resp. $\mathcal{UK}^0_{H}$) the class of all functions $f \in \mathcal{S}_{H}$ (resp. $f \in \mathcal{S}^0_{H}$) that are  uniformly convex in $\mathbb{D}$. Here a locally univalent harmonic function $f=h+\overline{g} \in \mathcal{H}$ is said to be uniformly convex in $\mathbb{D},$  if $f$ is fully convex in $\mathbb{D}$ and maps every circular arc $\gamma$ contained in $\mathbb{D},$ with center in $\ID$ onto a convex arc $f(\gamma)$ (see \cite{PonSaiPraj}). We refer to \cite{goodman1, goodman2, ronning} for discussion on this topic for the analytic case. A necessary and sufficient conditions for a harmonic mapping to be uniformly convex is (see \cite{PonSaiPraj})
\begin{equation}\label{unifconv}
{\rm Re\,} \{P(z, \zeta)\} \geq 0 ~\mbox{ for }(z, \zeta) \in \mathbb{D}\times\mathbb{D} ~\mbox{ such that }~ z \neq \zeta,
\end{equation}
where
$$P(z, \zeta) = \frac{(z-\zeta) h'(z)+(z-\zeta)^2 h''(z)+\overline{(z-\zeta)g'(z)}+ \overline{(z-\zeta)^2 g''(z)}}{(z-\zeta) h'(z)- \overline{(z-\zeta) g'(z)}}.
$$
Moreover, the following sufficient condition has been proved in \cite{PonSaiPraj}.

\begin{Thm}\label{jkps_T5}
Let $f =h+\overline{g}$ have the form \eqref{intro1} with $|b_1|<1$,  and satisfy the condition
\begin{equation}\label{unifconvsuf}
\sum_{n=2}^{\infty}n(2n-1)(|a_n|+|b_n|) \leq 1 -|b_1|.
\end{equation}
Then $f \in \mathcal{UK}_{H}$. The bound in \eqref{unifconvsuf} is sharp, especially for
$f(z)=z- \overline{e^{i3\alpha}z^2/6}$ in $\mathcal{UK}_{H}^0$, where $\alpha \in \IR$.
\end{Thm}

Using Theorem \Ref{jkps_T5} and  Gaussian hypergeometric functions, a family of examples of functions in $\mathcal{UK}^0_{H}$ is presented in Section \ref{p8sec5}.

It is easy to see that
$$ \mathcal{UK}_{H} \subset \widetilde{\mathcal{AK}_H} \subset \mathcal{AK}_H  \subset \mathcal{FK}_H.
$$
Further analysis shows that the classes $\widetilde{\mathcal{AK}_H}$ and $\mathcal{AK}_H$ are linear invariant families of univalent harmonic mappings, i.e., whenever $f = h +\overline{g} \in \mathcal{AK}_H$ (resp.  $\widetilde{\mathcal{AK}_H}$), the function
$$ F(z) = \frac{f(e^{i\theta}\frac{z+a}{1+\overline{a}z})-f(ae^{i\theta})}{(1-|a|^2)h'(ae^{i\theta})e^{i\theta}}
$$
also belongs to the class $\mathcal{AK}_H$ (resp. $\widetilde{\mathcal{AK}_H}$) for all $a \in \ID$ and $\theta \in \IR$. As disk automorphism takes circles onto circles, one could infer this fact from the definitions. These families of functions have the affine invariance property too. That is, if $f \in \mathcal{AK}_H$ (resp.  $\widetilde{\mathcal{AK}_H}$), then the function $(f + c \overline{f})/(1+cb_1) \in \mathcal{AK}_H$ (resp.  $\widetilde{\mathcal{AK}_H}$) for all $c \in \mathbb{D}$, where $b_1=g'(0)$. As a consequence of these properties, we derive many interesting results for the class $\mathcal{AK}_H$.

The article is organized as follows: In Section \ref{p8sec2}, we prove a necessary and sufficient condition for a sense-preserving harmonic function to be absolutely convex in $\mathbb{D}$. In Section \ref{p8sec3}, we find the order of the families $\mathcal{AK}_{H}$ (resp. $\mathcal{UK}_H$) and as an application, we prove growth and covering theorems for these families. In Section \ref{p8sec4}, we recall some results about uniformly starlike harmonic mappings and then discuss  Bernardi type integral transforms for harmonic functions. Using this we investigate its connection between a subclass of uniformly starlike harmonic functions and a subclass of uniformly convex harmonic functions. In Section \ref{p8sec5}, we present results leading to examples of uniformly starlike (resp. uniformly convex) harmonic mappings in $\ID$.

\section{Absolutely Convex Harmonic Mappings}\label{p8sec2}

In \cite{Duren_Osgood}, a necessary and sufficient condition for a harmonic function $f$ to belong to the class $\mathcal{FK}_{H}$ has been proved and it is as
follows:

\begin{Thm}\label{Thm_FK_NS}
Let $f=h+\overline{g}$ be a sense-preserving harmonic mapping of $\ID$. Then $f \in \mathcal{FK}_{H}$ if and only if
\beq\label{FK_NS}
|zh'(z)|^2{\rm Re\,}\left\{ 1+  \frac{zh''(z)}{h'(z)}\right\} &\geq& |zg'(z)|^2{\rm Re\,}\left\{ 1+  \frac{zg''(z)}{g'(z)}\right\} \\
\nonumber
&& \quad +{\rm Re\,}\left\{ z^3 [h''(z)g'(z)-h'(z)g''(z)]\right\}
\eeq
for all $z \in \ID$.
\end{Thm}

Now, we prove a necessary and sufficient condition for a sense-preserving harmonic function $f=h+\overline{g}$ to belong to the class $\mathcal{AK}_H$.

\bthm\label{p8thm1}
Let $f=h+\overline{g} \in \mathcal{H}$ be a sense-preserving harmonic function in $\ID$. Then the following statements are equivalent:
\begin{enumerate}
  \item[{\rm (i)}] $f \in \mathcal{AK}_H$
  \item[{\rm (ii)}] $f \in \widetilde{\mathcal{AK}_H}$
  \item[{\rm (iii)}]
  $\ds |h'(\zeta)|^2\left(1+|b|^2 + {\rm Re\,}\left\{ \frac{(\zeta - b)(1-\zeta\overline{b})h''(\zeta)}{h'(\zeta)} - 2 \overline{b}\zeta \right\} \right) \geq $\\
  $\ds |g'(\zeta)|^2 \left(1+|b|^2 + {\rm Re\,}\left\{ \frac{(\zeta - b)(1-\zeta\overline{b})g''(\zeta)}{g'(\zeta)} - 2 \overline{b}\zeta \right\} \right) \\
+ {\rm Re\,} \left\{ \frac{(\zeta - b)^3(1-\zeta \overline{b})^3}{|\zeta - b|^2 |1-\zeta \overline{b}|^2}
 \right\}$ for all $\zeta,$ $b \in \ID$.
\end{enumerate}
\ethm
\bpf
From the definition of $\mathcal{AK}_H$ and $\widetilde{\mathcal{AK}_H}$, it is clear that $\widetilde{\mathcal{AK}_H} \subset \mathcal{AK}_H$ and thus, we first prove that (i) $\Longrightarrow$ (ii).

Assume the contrary that  $f$ belongs to $\mathcal{AK}_H$ but not in $\widetilde{\mathcal{AK}_H}$. Then, there exists a disk $\ID(a, r) \subset \ID$ such that $|a|+r=1$ and $f(\ID(a, r))$ is not a convex domain. Therefore, there exist $w_1$, $w_2 \in f(\ID(a, r))$ such that
$$[w_1, w_2]:=\{w:\,w=tw_2+(1-t)w_1, ~0<t<1 \} \nsubseteq f(\ID(a, r)).
$$
Let $z_1 = f^{-1}(w_1)$ and $z_2 = f^{-1}(w_2)$. As $z_1,$ $z_2 \in \ID(a, r)$, we can find a $\rho > 0$ such that $\max\{|z_1 - a|, |z_2 - a| \} < \rho < r$. Since
$z_1,$ $z_2 \in \ID(a, \rho)$ and $f \in \mathcal{AK}_H$, it follows that $f(\ID(a, \rho))$ is convex, which is a contradiction to the fact that $[w_1, w_2] \nsubseteq f(\ID(a, r))$. Therefore, $f \in \widetilde{\mathcal{AK}_H}$. This completes the first part of the proof.

Next, we prove that (i) $\Longleftrightarrow$ (iii). We begin to let $f \in \mathcal{AK}_H$ and define $F$  by
$$\ds F(z) = \frac{f(e^{i\theta} \frac{z+a}{1+\overline{a}z})-f(ae^{i\theta})}{(1-|a|^2)e^{i\theta}h'(ae^{i\theta})} = H(z) + \overline{G(z)},
$$
where $H$ and $G$ are analytic in $\ID$ and
$$H'(z) = \frac{h'(e^{i\theta} \frac{z+a}{1+\overline{a}z})}{(1+\overline{a}z)^2h'(ae^{i\theta})} ~\mbox{ and }~ G'(z) = \frac{g'(e^{i\theta} \frac{z+a}{1+\overline{a}z})e^{2i\theta}}{(1+\overline{a}z)^2\overline{h'(ae^{i\theta})}}.
$$
From the definition of the class $\mathcal{AK}_H$, we have $f \in \mathcal{AK}_H$ if and only if $F \in \mathcal{FK}_H$ for all $a \in \ID$ and $\theta \in \IR$.
From equation \eqref{FK_NS} of Theorem \Ref{Thm_FK_NS}, it is clear that $f \in \mathcal{AK}_H$ if and only if
\beq\label{AK_NS}
|zH'(z)|^2{\rm Re\,}\left\{ 1+  \frac{zH''(z)}{H'(z)}\right\} &\geq& |zG'(z)|^2{\rm Re\,}\left\{ 1+  \frac{zG''(z)}{G'(z)}\right\} \\ \nonumber
&& \quad + {\rm Re\,}\left\{ z^3 [H''(z)G'(z)-H'(z)G''(z)]\right\}.
\eeq
As
$$H''(z) = \frac{h''(e^{i\theta} \frac{z+a}{1+\overline{a}z})(1-|a|^2)e^{i\theta} - h'(e^{i\theta} \frac{z+a}{1+\overline{a}z})2\overline{a}(1+\overline{a}z) }{(1+\overline{a}z)^4h'(ae^{i\theta})}
$$
and
$$G''(z) = \frac{g''(e^{i\theta} \frac{z+a}{1+\overline{a}z})(1-|a|^2)e^{i\theta} - g'(e^{i\theta} \frac{z+a}{1+\overline{a}z})2\overline{a}(1+\overline{a}z) }{(1+\overline{a}z)^4\overline{h'(ae^{i\theta})}}e^{2i\theta},
$$
the inequality \eqref{AK_NS} is equivalent to
\beqq
\left|\frac{\zeta - b}{1-\zeta\overline{b}} \right|^2 \left|\frac{h'(\zeta)(1-\zeta\overline{b})^2}{h'(b)} \right|^2 {\rm Re\,}\left\{1 + \left(\frac{\zeta - b}{1-\zeta\overline{b}} \right) \left( \frac{h''(\zeta)(1-\zeta \overline{b})^2 - h'(\zeta)2\overline{b}(1-\zeta \overline{b})}{h'(\zeta)(1-|b|^2)} \right)\right\} \geq \\
\left|\frac{\zeta - b}{1-\zeta\overline{b}} \right|^2 \left|\frac{g'(\zeta)(1-\zeta\overline{b})^2}{h'(b)} \right|^2 {\rm Re\,}\left\{1 + \left(\frac{\zeta - b}{1-\zeta\overline{b}} \right) \left( \frac{g''(\zeta)(1-\zeta \overline{b})^2 - g'(\zeta)2\overline{b}(1-\zeta \overline{b})}{g'(\zeta)(1-|b|^2)} \right)\right\} \\
+ {\rm Re\,}\left\{(\zeta - b)^3 \left[ \frac{\{h''(\zeta)g'(\zeta) - g''(\zeta)h'(\zeta)\}(1-\zeta \overline{b})^3}{|h'(b)|^2(1-|b|^2)} \right]\right\} ~\mbox{ for all }~ \zeta, b \in \ID,
\eeqq
where $\zeta = e^{i\theta}(z+a)/(1+\overline{a}z)$  and $b = a e^{i\theta}.$
On simplification, the above inequality reduces to
\beqq
|h'(\zeta)(1-\zeta \overline{b})(\zeta - b)|^2 \left(1+|b|^2 + {\rm Re\,}\left\{ \frac{(\zeta - b)(1-\zeta\overline{b})h''(\zeta)}{h'(\zeta)} - 2 \overline{b}\zeta \right\} \right) \geq \\
|g'(\zeta)(1-\zeta \overline{b})(\zeta - b)|^2 \left(1+|b|^2 + {\rm Re\,}\left\{ \frac{(\zeta - b)(1-\zeta\overline{b})g''(\zeta)}{g'(\zeta)} - 2 \overline{b}\zeta \right\} \right) \\
+ {\rm Re\,} \left\{ (\zeta - b)^3(1-\zeta \overline{b})^3[h''(\zeta)g'(\zeta)-h'(\zeta)g''(\zeta)]\right\},
\eeqq
which can be rewritten in the form of inequality (iii) in the statement.
\epf

\section{Linear and Affine Invariant Families}\label{p8sec3}
In \cite{Starkov_2004}, the author gave a new definition for the order of a linear invariant family $L$ of harmonic mappings
$f$ have the form \eqref{intro1} as
$$ \overline{{\rm ord}}~L = \sup_{f\in L} \frac{|a_2-\overline{b_1}b_2|}{1-|b_1|^2}.
$$
and proved upper and lower bounds for the Jacobian of $f$ and many more interesting results (see also \cite{Starkov_2011}).

\begin{Lem}{\rm (cf. \cite{Starkov_Ganenkova})}\label{new_ord}
Let L be a linear invariant family of harmonic mappings. Then,
$\overline{{\rm ord}}~L = \sup_{f\in A^0[L]} |a_2|,$
where
$$A^0[L] = \left\{F=\frac{f+\epsilon \overline{f}}{1+\epsilon \overline{b_1}}:\, f\in L, \epsilon \in \ID, F_{\overline{z}}(0)=0  \right\}.
$$
\end{Lem}

\begin{Thm}{\rm (cf. \cite{Starkov_2011})}\label{Thm_Jacob1}
Let $f \in L$ with $b_1=f_{\overline{z}}(0)$, where $L$ is as above. Then, the Jacobian $J_f$ of the mapping $f$ with any $z \in \mathbb{D}$ satisfies the
bounds
$$
(1-|b_1|^2)\frac{(1-|z|)^{2\alpha_0 - 2}}{(1+|z|)^{2\alpha_0 + 2}} \leq J_f(z) \leq
(1-|b_1|^2)\frac{(1+|z|)^{2\alpha_0 - 2}}{(1-|z|)^{2\alpha_0 + 2}},
$$
where $\alpha_0 = \overline{{\rm ord}}~L$.
\end{Thm}

Next, we shall use the affine and linear invariance properties of $\mathcal{AK}_{H}$ to prove growth, covering and area
theorems.

\bthm\label{p8thm2}
Suppose that $f$ is of the form \eqref{intro1} and $f \in \mathcal{AK}^{0}_{H}$. Then the coefficient $a_2$ satisfies the condition $|a_2|  \leq 2/\sqrt{3} \approx 1.1547 $.
Every function
$f \in \mathcal{AK}^{0}_{H}$ satisfies the inequalities
\be\label{Growth_Ineq_AK}
\left(\frac{\sqrt{3}}{\sqrt{3}+4}\right)\left[1-\left(\frac{1-r}{1+r}\right)^{\frac{\sqrt{3}+4}{2\sqrt{3}}} \right] \leq |f(z)| \leq
\left(\frac{\sqrt{3}}{\sqrt{3}+4}\right) \left[\left(\frac{1+r}{1-r}\right)^{\frac{\sqrt{3}+4}{2\sqrt{3}}} -1 \right],
\ee
for $0 < |z|= r < 1$. In particular, the range of each function $f \in \mathcal{AK}^{0}_{H}$ contains the disk $|w| < \sqrt{3}/(\sqrt{3}+4) \approx 0.302169$.
\ethm
\bpf
Let $f = h +\overline{g} \in \mathcal{AK}^{0}_{H}$ be of the form \eqref{intro1}. Without loss of generality, we may assume that $a_2 > 0$.
In Theorem \ref{p8thm1} (iii), we may let $\zeta = 0$, and make use of the
facts that $h'(0)=1$ and $g'(0)=0$. This gives the inequality
$$ {\rm Re\,} \{a_2 e^{i t}+b_2e^{3i t}\} \leq \frac{1+|b|^2}{2|b|} ~\mbox{ for all }~t\in\IR,$$
where $b = |b|e^{i t} \neq 0$. Since the above inequality holds for each $b \in \ID$, we have
$$ {\rm Re\,} \{a_2 e^{i t}+b_2e^{3i t}\} \leq 1 ~\mbox{ for all }~t\in\IR .$$
Set $t_k=t+(k/3)\pi$, for $t\in\IR$ and $k=0,1$. From the above inequality, it is clear that
$$ {\rm Re\,} \{a_2 e^{i t_k}+b_2e^{3i t_k}\} \leq 1 ~\mbox{ for all }~t_k\in\IR ~\mbox{ and }~k=0,1,$$
which shows that
$$ {\rm Re\,} \{a_2 e^{i t}+b_2e^{3i t}\} + {\rm Re\,} \{a_2 e^{i (t+\pi/3)}- b_2e^{3i t}\} = {\rm Re\,} \{a_2 e^{it}(1+e^{i\pi/3})\} \leq 2
$$
for all $t\in\IR$. Now, we choose $t$ such that $a_2 e^{it}=|a_2| e^{-i\pi/6}$. Then the last inequality is equivalent to $|a_2| \leq 2/\sqrt{3}$.
Since every function in the class $\mathcal{AK}_{H}$ can be obtained by composing a function from the class $\mathcal{AK}^{0}_{H}$ with affine mapping,
it is easy to see that for functions  $f = h +\overline{g} \in \mathcal{AK}_{H}$ with the representation as in \eqref{intro1}, the corresponding second coefficient $a_2$ satisfies the inequality
$$|a_2| \leq \frac{2}{\sqrt{3}} + \frac{1}{2} = \frac{\sqrt{3}+4}{2\sqrt{3}}.$$
From a result of Sheil-Small \cite{Sheil-Small} on linear invariant families of univalent harmonic mappings, we obtain the following inequality
$$
\left(\frac{\sqrt{3}}{\sqrt{3}+4}\right)\left[1-\left(\frac{1-r}{1+r}\right)^{\frac{\sqrt{3}+4}{2\sqrt{3}}} \right] \leq |f(z)| \leq
\left(\frac{\sqrt{3}}{\sqrt{3}+4}\right) \left[\left(\frac{1+r}{1-r}\right)^{\frac{\sqrt{3}+4}{2\sqrt{3}}} -1 \right],
$$
whenever $f \in \mathcal{AK}^{0}_{H}$.
Allowing $r \rightarrow 1^-$ in the left hand side of above inequality, we conclude that the range of $f$ contains a disk of radius at least $\sqrt{3}/(\sqrt{3}+4)$.
\epf

\br
If $f$ is of the form \eqref{intro1} and $f \in \mathcal{AK}^{0}_{H}$ with $a_2 \geq 0$, then $a_2 + {\rm Re\,}b_2 \leq 1$. This inequality is
sharp in the case of analytic functions.
\er

\bthm\label{p8thm2.1}
Let $f \in \mathcal{AK}_{H}$ with $b_1=f_{\overline{z}}(0)$. Then  the Jacobian $J_f$ of the mapping $f$ with any $z \in \mathbb{D}$ satisfies the
bounds
\be\label{Jacob_Ineq1}
(1-|b_1|^2)\frac{(1-|z|)^{\frac{4}{\sqrt{3}} - 2}}{(1+|z|)^{\frac{4}{\sqrt{3}} + 2}} \leq J_f(z) \leq
(1-|b_1|^2)\frac{(1+|z|)^{\frac{4}{\sqrt{3}} - 2}}{(1-|z|)^{\frac{4}{\sqrt{3}} + 2}}.
\ee
Moreover, $A(f(\mathbb{D}_r))$, the area of $f(\mathbb{D}_r)$ have the following bounds.
$$
\frac{\pi (1-|b_1|^2)}{26}\left\{3 - (3r^2 + 8\sqrt{3}r +3)\frac{(1-r)^{\frac{4}{\sqrt{3}}-1}}{(1+r)^{\frac{4}{\sqrt{3}}+1}}  \right\} \leq A(f(\mathbb{D}_r)) $$
and
$$
A(f(\mathbb{D}_r)) \leq \frac{\pi (1-|b_1|^2)}{26}\left\{3 - (3r^2 - 8\sqrt{3}r +3)\frac{(1+r)^{\frac{4}{\sqrt{3}}-1}}{(1-r)^{\frac{4}{\sqrt{3}}+1}}  \right\}.
$$
\ethm
\bpf
From Lemma \Ref{new_ord} and Theorem \ref{p8thm2}, it follows that $\overline{{\rm ord}}~\mathcal{AK}_{H} = 2/\sqrt{3}$. The proof for the inequalities \eqref{Jacob_Ineq1} follows from Theorem \Ref{Thm_Jacob1}. To find the upper bound and lower bound for $A(f(\mathbb{D}_r))$, let us consider
\beqq
A(f(\mathbb{D}_r)) &=& \iint\limits_{\mathbb{D}_r} \,J_f(z) dx\,dy\\
&\leq& (1-|b_1|^2) \int_0^{2\pi} \!\!\!\int_0^r \frac{(1+\rho)^{\frac{4}{\sqrt{3}} - 2}}{(1-\rho)^{\frac{4}{\sqrt{3}} + 2}} \rho \, d\rho\,d\theta\\
&=& 2\pi (1-|b_1|^2) \int_0^r \frac{(1+\rho)^{\frac{4}{\sqrt{3}} - 2}}{(1-\rho)^{\frac{4}{\sqrt{3}} + 2}} \rho \, d\rho\\
&=& \frac{\pi (1-|b_1|^2)}{26}\left\{3 - (3r^2 - 8\sqrt{3}r +3)\frac{(1+r)^{\frac{4}{\sqrt{3}}-1}}{(1-r)^{\frac{4}{\sqrt{3}}+1}}  \right\}.
\eeqq
The lower bound for $A(f(\mathbb{D}_r))$ follows in a similar way.
\epf

\bthm\label{p8thm3}
Suppose that $f$ is of the form \eqref{intro1} and $f \in \mathcal{UK}^{0}_{H}$. Then the coefficient $a_2$ satisfies the condition $|a_2|  \leq 1/\sqrt{3} \approx 0.57735 $.
Every function
$f \in \mathcal{UK}^{0}_{H}$ satisfies the inequalities \eqref{Growth_Ineq_AK}.
\ethm
\bpf
Let $f$ be of the form \eqref{intro1} and $f \in \mathcal{UK}^{0}_{H}$. By setting $z=0$ and allowing $|\zeta| \rightarrow 1^{-}$ in \eqref{unifconv}, we have
$$ {\rm Re\,} \{a_2 e^{i t}+\overline{b_2}e^{-3i t}\} \leq \frac{1}{2} ~\mbox{ for all }~t\in\IR .$$
The rest of the proof is similar to the proof of Theorem \ref{p8thm2} and hence we omit the details.
\epf

\br
The coefficient bound obtained in Theorem \ref{p8thm3} serves as a necessary condition to check whether a function belongs to the class
$\mathcal{UK}^{0}_{H}$. For a harmonic function $f \in \mathcal{UK}^{0}_{H}$, we can not improve the bounds in the inequalities \eqref{Growth_Ineq_AK},
as $\mathcal{UK}_{H}$ is not a linear invariant family.
\er

\section{Alexander Type Theorems for the subclasses $\mathcal{US}^*_{H}$ and $\mathcal{UK}_{H}$}\label{p8sec4}

A complex valued function $f \in C^1(\mathbb{D})$ is said to be uniformly starlike in $\mathbb{D}$ (see \cite{PonSaiPraj}) if $f$ is fully starlike in $\mathbb{D}$ and maps every circular arc $\gamma_{\zeta}$ contained in $\mathbb{D}$ with center $\zeta \in \mathbb{D},$ to the arc $f(\gamma_{\zeta})$ which is starlike with respect to $f(\zeta)$.

Let $\mathcal{US}^*_{H}$ (resp. $\mathcal{US}^{0*}_{H}$) denote the class of all mappings $f \in \mathcal{S}_{H}$ (resp. $ f \in \mathcal{S}^0_{H}$) that are  uniformly starlike in $\mathbb{D}$. Recently,  in \cite{PonSaiPraj} these classes have been introduced and studied. In particular, the following results are derived.

\begin{Thm}\label{jkps_T1}{\rm \cite{PonSaiPraj}}
Let $f \in C^1(\mathbb{D})$ be such that
\begin{itemize}
  \item [{\rm (i)}]$f(0)=0$ and $J_{f}(z) > 0$ for all $z \in \mathbb{D}$
  \item [{\rm (ii)}] ${\rm Re} \displaystyle \left\{\frac{f(z) - f(\zeta)}{(z-\zeta)f_z(z) - \overline{(z-\zeta)}f_{\overline{z}}(z)}\right\} \geq 0 ~\mbox{ for }~ (z, \zeta) \in \mathbb{D}\times\mathbb{D} ~\mbox{ such that }~ z \neq \zeta$ and the strict inequality holds when $\zeta=0.$
\end{itemize}
Then $f$ is univalent and uniformly starlike in $\mathbb{D}$. In particular, if $f \in {\mathcal{H}}$ and satisfies the conditions {\rm (i)} and {\rm (ii)}, then $f \in \mathcal{US}^*_{H}$.
\end{Thm}

\begin{Thm}{\rm \cite{PonSaiPraj}}\label{jkps_T3}
Let $f =h+\overline{g}$ have the form \eqref{intro1} and satisfy the condition
\begin{equation}\label{th2.1}
\sum_{n=2}^{\infty} n |a_n| + \sum_{n=1}^{\infty} n |b_n| \leq \frac{1}{2}.
\end{equation}
Then $f \in \mathcal{US}_{H}^*.$
\end{Thm}

Alexander's theorem for analytic function states that $f \in \mathcal{K} \Longleftrightarrow z f' \in \mathcal{S}^*$, but
there does not exist a similar two way implications between the classes
$$UCV := \mathcal{UK}_{H}\cap\{f =h+\overline{g}:\, g \equiv 0 \}  ~\mbox{ and  }~UST := \mathcal{US}^*_{H} \cap \{f =h+\overline{g}:\, g \equiv 0 \}
$$
(see \cite{goodman2, ronning}). It is possible to present a one way bridge between the subclasses of $\mathcal{US}^*_{H}$ and $\mathcal{UK}_{H}$. Now, we introduce an operator $f=h+\overline{g} \mapsto \Lambda_f = \Lambda_h + \overline{\Lambda_g}$ taking function $f\in \mathcal{US}^*_{H}$ into $\Lambda_f \in \mathcal{UK}_{H}$. Our particular emphasize will be when $\Lambda_h = H_{a, b}(h)$ for $a, b \in \IR$, $a \ne b$, $a > -1$ and $b>-1$, where
$$
H_{a, b}(h)(z)=\frac{(a+1)(b+1)}{(b-a)}\int\limits_0^1 t^{a - 1} (1-t^{b-a}) h(tz) \, \mathrm{d}t = \sum_{n=1}^{\infty} \frac{(a+1)(b+1)}{(a+n)(b+n)}a_n z^n,
$$
where $h(z)=z+\sum_{n=2}^{\infty}a_n z^n$.

These integral transforms actually leads to the convolution of $f$ with certain classes of special functions (see \cite{PonRon_1997}). In particular, we get Bernardi type operator in the limiting case when $b \rightarrow \infty$. As in \cite{PonRon_1997}, here are two limiting cases of $H_{a, b}(h)$ given by
\be\label{p8aeq1}
H_{a, \infty}(h)(z)= \lim_{b \to \infty} H_{a, b}(h)(z) = \frac{a + 1}{z^{a}} \int\limits_0^z t^{a - 1} h(t) \, \mathrm{d}t ,
\ee
and
\be\label{p8aeq2}
H_{a, a}(h)(z)= \lim_{b \to a} H_{a, b}(h)(z) = -(a+1)^2\int\limits_0^1 t^{a - 1} (\log t)\, h(tz) \, \mathrm{d}t .
\ee
Note that
$$ H_{a, \infty}(h)(z) = \sum_{n=1}^{\infty} \frac{a+1}{a+n}a_n z^n ~\mbox{ and }~  H_{a, a}(h)(z) = \sum_{n=1}^{\infty} \frac{(a+1)^2}{(a+n)^2}a_n z^n.
$$

\bthm\label{jkps_T6}
Suppose that $f =h+\overline{g} \in \mathcal{H}$ is of the form \eqref{intro1}
and satisfies the coefficient condition \eqref{th2.1}. Set $H = \Lambda_h$ and $G=\Lambda_g$.
The harmonic function $F(z)=H(z)+\overline{G(z)} \in \mathcal{UK}_{H}$, provided $a>-1$, $b>-1$ and
satisfy one of the following conditions:
\begin{enumerate}
  \item[{\rm (i)}] $ab \le 3$.\\
  \item[{\rm (ii)}] $ab > 3$ and $a^2b^2-4ab-2(a+b) \le 1$.
\end{enumerate}
\ethm
\bpf
Suppose that $f$ is of the form \eqref{intro1} and satisfies the inequality \eqref{th2.1}, i.e.,
$$
\sum_{n=2}^{\infty} 2n (|a_n| + |b_n|) \leq 1 - 2|b_1|.
$$
We may write $F(z)$ as
$$F(z) = \sum_{n=1}^{\infty} A_{n} z^n +  \overline{\sum_{n=1}^{\infty} B_{n} z^n},$$ where
$$
A_{n} = \frac{(a+1)(b+1)}{(a+n)(b+n)} a_n ~\mbox{ and }~ B_{n} = \frac{(a+1)(b+1)}{(a+n)(b+n)} b_n ~\mbox{ for }~ n \geq 1.
$$
From \eqref{unifconvsuf}, it suffices to show that the coefficients of $F(z)$ satisfy the condition
\begin{equation}\label{unifconvsuf1}
\sum_{n=2}^{\infty}n(2n-1)(|A_n|+|B_n|) \leq 1 -|B_1|.
\end{equation}
We shall show that this holds provided either (i) or (ii) hold. It is a simple exercise to observe that \eqref{unifconvsuf1} holds whenever $\varphi(n) \ge 0$
for all $n \ge 2$, where
\beqq\nonumber
\varphi(n)&=&2n^2-2(1+ab)n+(3ab+a+b+1)\\
&=& 2(n-2)^2 + 2(n-2)(3-ab)+(a+b+2)+(3-ab).
\eeqq
If $ a, b >-1$ and $3 -ab \ge 0$, then $\varphi(n) > 0$ for all $n \ge 2$.

Next, we consider the case $3 -ab < 0$. In this case, we rewrite
$$\varphi(n) = 2\left\{\left[(n-2)-\frac{ab-3}{2}\right]^2 - \left[\frac{(ab-3)^2-2(a+b+5-ab)}{4}\right] \right\}
$$
which gives that $\varphi(n) \ge 0$ for all $n \ge 2$, whenever
$$(ab-3)^2-2(a+b+5-ab) \le 0.$$
Thus, under (ii), the desired conclusions follow.
\epf

\br
In the proof of Theorem \ref{jkps_T6}, if ~$a b>3$, $a^2b^2-4ab-2(a+b) > 1$ and $\lfloor r_1 \rfloor - \lfloor r_2 \rfloor =0$,
then the function $F(z)$ defined in Theorem \ref{jkps_T6} belongs to the class $\mathcal{UK}_{H}$, where $r_1$ and $r_2$
are real roots of the equation $\varphi(n)=0$ and $\lfloor x \rfloor$ denotes the greatest integer less than or equal to $x$.
For example, if we take $a=2$ and $b=59/20$,
then neither the condition (i) nor (ii) of Theorem \ref{jkps_T6} is satisfied. But, in this case
$$ \varphi(n) = 2 n^2 - \frac{69}{5}n + \frac{473}{20}  >  0 ~\mbox{ for all }~ n \ge 2.
$$
Consequently, the corresponding $F(z)$ with $a=2$ and $b=59/20$ is in $\mathcal{UK}_{H}$.
\er

We now consider two special cases of Theorem \ref{jkps_T6}.\\

\noindent{\bf Case (i)} If we allow $b \rightarrow \infty$ in Theorem \ref{jkps_T6}, then the harmonic function $F(z)$ reduces to $F(z)=H_{a, \infty}(h)(z)+\overline{H_{a, \infty}(g)(z)}$, where $H_{a, \infty}(h)$ is defined by \eqref{p8aeq1}.
We conclude that the harmonic function $F(z)$ is univalent and uniformly convex in $\ID$, whenever $a \in (-1, ~0]$.\\

\noindent{\bf Case (ii)} If we take $a = b$ in Theorem \ref{jkps_T6}, then the harmonic function $F(z)$ reduces to
$F(z)=H_{a, a}(h)(z)+\overline{H_{a, a}(g)(z)}$, where $H_{a, a}(h)$ is defined by \eqref{p8aeq2}.
It follows that the harmonic function $F(z)$ is univalent and uniformly convex in $\ID$, whenever
$$a \ds \in \left(-1, ~\frac{3\sqrt{2}-\sqrt{5}}{\sqrt{5}-\sqrt{2}}\right].
$$

Indeed we have $F(z) = z +\sum_{n=2}^{\infty} A_{n} z^n + \overline{\sum_{n=1}^{\infty} B_{n} z^n}$, where
\be\label{3.5}
A_{n} = \left(\frac{1+a}{n+a}\right)^2 a_n ~\mbox{ and }~ B_{n} = \left(\frac{1+a}{n+a}\right)^2 b_n ~\mbox{ for }~ n \geq 1.
\ee
As in the proof of Theorem \ref{jkps_T6}, for $F$ to be univalent and uniformly convex in $\ID$, it suffices to verify the inequality \eqref{unifconvsuf},
with $A_n$ and $B_n$ as given in \eqref{3.5}, which indeed holds provided
$$(a+1)^2(2n-1) \le 2(a+n)^2 ~\mbox{ for all }~ n \ge 2,
$$
i.e., if
$$ a \le \frac{\sqrt{2} n - \sqrt{2n-1}}{\sqrt{2n-1} - \sqrt{2}}=:\psi(n) ~\mbox{ for all }~ n \ge 2.
$$
It is easy to see that
$$
\min_{n \in \IN \setminus \{1\}} \psi(n) =\psi(3) = \frac{3\sqrt{2}-\sqrt{5}}{\sqrt{5}-\sqrt{2}},
$$
which gives the stated range for $a$. This is handy instead of considering the conditions (i)--(ii) of Theorem \ref{jkps_T6} with $a=b$.

\bthm\label{jkps_T7}
Suppose that $F =H+\overline{G} \in \mathcal{H}$, $G'(0)=0$ and satisfies the coefficient condition \eqref{unifconvsuf}. For $h_{a}$ and $g_{a}$ defined by
$$h_{a}(z) = \frac{a H(z)+zH'(z)}{a+1} ~\mbox{ and }~ g_{a}(z) = \frac{a G(z)+zG'(z)}{a+1},$$
where $a \geq 1$, the harmonic function $f_{a}(z)=h_{a}(z)+\overline{g_{a}(z)} \in \mathcal{US}^*_{H}$.

\ethm
\bpf
The proof is similar to the proof of Theorem \ref{jkps_T6} and so we omit its details.
\epf

\section{Applications and Examples}\label{p8sec5}
Here, we consider harmonic mappings whose co-analytic parts involve Gaussian hypergeometric function. The Gaussian hypergeometric function $ _{2}F_{1}(a,b;c;z) $ is defined as
\begin{eqnarray}\label{hyper1}
_{2}F_{1}(a,b;c;z) := F(a,b;c;z)=\sum_{n=0}^{\infty} \frac{(a)_n (b)_n}{(c)_n \,n!}z^n, ~|z| < 1,
\end{eqnarray}
where $a,b,c \in \mathbb{C}$, $c \neq 0,-1,-2, \cdots ,$ and $(a)_n$ is the Pochhammer symbol defined by $\,(a)_0=1$ \,and $\,(a)_n=a(a+1)\cdots(a+n-1) \;(n \in \mathbb{N}).$ The series \eqref{hyper1} is absolutely convergent for $|z|<1$. Moreover, if ${\rm Re\,}(c)>{\rm Re\,}(a+b)$, then it is also convervent for $|z| \leq 1.$ The following well-known Gauss formula \cite{temme} and Lemma \Ref{hyper_form} given below are crucial in the proof of our results of this section:
\begin{eqnarray}\label{hyper2}
F(a,b;c;1)= \frac{\Gamma(c) \Gamma(c-a-b)}{\Gamma(c-a) \Gamma(c-b)} < \infty \;\;\;\mbox{for} \;\;{\rm Re\,}(c)>{\rm Re\,}(a+b).
\end{eqnarray}

\begin{Lem}\label{hyper_form}{\rm \cite[Lemma 3.1]{PonRon} and \cite{kim_samy}}. Let $a, b >0.$ Then we have the following
\begin{enumerate}
\item[(a)] For $\,c>a+b+1,$
$$\sum_{n=0}^{\infty} \frac{(n+1)(a)_n (b)_n}{(c)_n \,n!} = \frac{\Gamma(c) \Gamma(c-a-b-1)}{\Gamma(c-a)\Gamma(c-b)}(ab+c-a-b-1).$$
\item[(b)] For $\,c>a+b+2,$
$$\sum_{n=0}^{\infty} \frac{(n+1)^2(a)_n (b)_n}{(c)_n \,n!} = \frac{\Gamma(c) \Gamma(c-a-b)}{\Gamma(c-a)\Gamma(c-b)}\left[\frac{(a)_2(b)_2}{(c-a-b-2)_2}+\frac{3ab}{c-a-b-1}+1 \right].$$
\end{enumerate}
\end{Lem}

\bthm\label{jkps_T8}
Let either $\,a,b \in (-1, \infty)$ with $ab>0,$ or $a,b \in \mathbb{C}\setminus \{0\}$ with $b=\overline{a}.$ Assume that $c$ is a positive real number, $\alpha \in \mathbb{D}$,
$$f_1(z)=z+\overline{(\alpha/2) z F(a,b;c;z)}, \;\;f_2(z)=z+\overline{(\alpha/2) (F(a,b;c;z)-1)} \quad \mbox{and}$$
$$f_3(z)=z+\overline{(\alpha/2) \int_0^z F(a,b;c;t) dt}. \qquad \qquad \qquad \qquad \qquad \qquad  \qquad \qquad $$
\begin{enumerate}
\item[{\rm (a)}] If $c > {\rm Re\,}(a+b)+1$ and
\begin{eqnarray}\label{hyper3}
|\alpha| \frac{\Gamma(c) \Gamma(c-a-b-1)}{\Gamma(c-a)\Gamma(c-b)}(ab+c-a-b-1) \leq 1,
\end{eqnarray}
then $f_1 \in \mathcal{US}_H^*.$
\item[{\rm (b)}]If $c > {\rm Re\,}(a+b)+1$ and
\begin{eqnarray}\label{hyper4}
ab|\alpha|\frac{\Gamma(c) \Gamma(c-a-b-1)}{\Gamma(c-a)\Gamma(c-b)}\leq 1,
\end{eqnarray}
then $f_2 \in \mathcal{US}_H^*.$
\item[{\rm (c)}] If $c > {\rm Re\,}(a+b)$ and
\begin{eqnarray}\label{hyper5}
|\alpha|\frac{\Gamma(c)\Gamma(c-a-b)}{\Gamma(c-a)\Gamma(c-b)}\leq 1,
\end{eqnarray}
then $f_3 \in \mathcal{US}_H^*.$
\end{enumerate}
\ethm
\begin{proof}

(a) Set $f_1(z)=z+\overline{g_1(z)}$, where $g_1(z)= (\alpha/2) z F(a,b;c;z)=\sum_{n=1}^{\infty} b_n z^n$.
We see that
$$
b_n = \frac{\alpha}{2} \frac{(a)_{n-1} (b)_{n-1}}{(c)_{n-1} \,(n-1)!} \;\;\mbox{for}\;\;n \geq 1,
$$
and
$$K:= \sum_{n=0}^{\infty} (n+1) |b_{n+1}| =\frac{|\alpha|}{2} \sum_{n=0}^{\infty} (n+1)\frac{(a)_{n} (b)_{n}}{(c)_{n} n!}.$$
From the formula (a) of Lemma \Ref{hyper_form}, \eqref{hyper3} is equivalent to $2 K\leq 1$ and thus $f_1 \in \mathcal{US}_H^*$,
by Theorem \Ref{jkps_T3}.\\

(b) For the proof of Case (b), we consider $g_2(z)=(\alpha/2)(F(a,b;c;z)-1)$
so that $f_2(z)=z+\overline{g_2(z)}$. Using \eqref{hyper2} it follows that
\begin{eqnarray*}
K &=& \frac{|\alpha|}{2} \sum_{n=1}^{\infty} n\frac{(a)_{n} (b)_{n}}{(c)_{n} \,n!} = \frac{|\alpha|}{2} \frac{ab}{c} \sum_{n=0}^{\infty} \frac{(a+1)_{n} (b+1)_{n}}{(c+1)_{n} n!}\notag \\
&=& \frac{|\alpha|ab}{2}\frac{\Gamma(c)\Gamma(c-a-b-1)}{\Gamma(c-a)\Gamma(c-b)} \leq \frac{1}{2}, \notag
\end{eqnarray*}
which is equivalent to \eqref{hyper4} and hence by the coefficient inequality \eqref{th2.1}, $f_2 \in \mathcal{US}_H^*.$\\

(c) For $f_3(z)$, we have $f_3(z)= z+\overline{g_3(z)}$, where
$$g_3(z) = \frac{\alpha}{2}\sum_{n=0}^{\infty} \frac{(a)_{n} (b)_{n}}{(c)_{n} \,n!}\frac{z^{n+1}}{(n+1)},$$
and thus, we find that
\begin{eqnarray*}
K &=&\frac{|\alpha|}{2}\sum_{n=0}^{\infty} \frac{(a)_{n} (b)_{n}}{(c)_{n} n!}=\frac{|\alpha|}{2} \frac{\Gamma(c)\Gamma(c-a-b)}{\Gamma(c-a)\Gamma(c-b)} \leq \frac{1}{2},\notag \\
\end{eqnarray*}
which is equivalent to \eqref{hyper5}. Again by Theorem \Ref{jkps_T3}, we conclude that $f_3 \in \mathcal{US}_H^*.$
\end{proof}

The case $a=1$ of Theorem \ref{jkps_T8} gives

\bcor Let $b$ and $c$ be positive real numbers and $\alpha \in \mathbb{D}$.
\begin{enumerate}
\item[{\rm (a)}] If
\begin{eqnarray}\label{hyper8}
c \geq \beta^+=\frac{2b+3-3|\alpha|+\sqrt{|\alpha|^2+2|\alpha|(2b^2-1)+1}}{2(1-|\alpha|)},
\end{eqnarray}
then the function $f(z)= z+ \overline{(\alpha/2) z F(1,b;c;z)}$ is uniformly starlike in $\mathbb{D}.$
\item[{\rm (b)}] If
\begin{eqnarray}\label{hyper9}
c \geq \gamma^+=\frac{2b+3+b|\alpha|+\sqrt{b^2(|\alpha|^2+4|\alpha|)+2|\alpha|b+1}}{2},
\end{eqnarray}
then the function $f(z)= z+ \overline{(\alpha/2) (F(1,b;c;z)-1)}$ is uniformly starlike in $\mathbb{D}.$
\item[{\rm (c)}] If $ c \geq 1+b/(1-|\alpha|),$
then function $f$ defined by
$$f(z)= z+ \overline{\frac{\alpha}{2} \int_0^z F(1,b;c;t)\, dt}$$ is uniformly starlike in $\mathbb{D}.$
\end{enumerate}
\ecor
\begin{proof}
(a) It suffices to prove that if $c \geq \beta^+,$  then the inequality \eqref{hyper3} is satisfied with $a=1.$ It can be easily seen that the condition $c \geq \beta^+$ implies that $c>b+2$. Next, the condition \eqref{hyper3} for $a=1$ reduces to
$$(1-|\alpha|)c^2+c(3|\alpha|-2b-3)+b^2+3b+2-2|\alpha| \geq 0.$$
Simplifying this inequality gives $(1-|\alpha|)(c-\beta^-)(c-\beta^+) \geq 0,$
where $\beta^+$ is given by \eqref{hyper8} and
\begin{eqnarray*}
\beta^-=\frac{2b+3-3|\alpha|-\sqrt{|\alpha|^2+2|\alpha|(2b^2-1)+1}}{2(1-|\alpha|)}.
\end{eqnarray*}
Since $\beta^+ \geq \beta^-$ and by hypothesis $c \geq \beta^+,$ the inequality \eqref{hyper3} holds and hence
$f(z)= z+ \overline{(\alpha/2) z F(1,b;c;z)} \in \mathcal{US}_H^*.$\\

(b) If we set $a=1$ in Theorem \ref{jkps_T8}(b), then the inequality \eqref{hyper4} is seen to be equivalent to
$$c^2-c(2b+3+b|\alpha|)+2 + 3 b + b |\alpha|+ b^2=(c-\gamma^-)(c-\gamma^+) \geq 0,$$
where $\gamma^+$ is given by \eqref{hyper9} and
\begin{eqnarray*}
\gamma^-=\frac{2b+3+b|\alpha|-\sqrt{b^2(|\alpha|^2+4|\alpha|)+2|\alpha|b+1}}{2}.
\end{eqnarray*}
Since $\gamma^+ \geq \gamma^-,$ the hypothesis that $c \geq \gamma^+$ gives the desired conclusion.\\

The proof of case (c) follows from Theorem \ref{jkps_T8}(c), if one adopts a similar approach.
\end{proof}
\noindent Note that for $\eta \in \mathbb{C} \setminus \{-1,-2,\cdots \}$ and $n \in \mathbb{N}_0=\mathbb{N} \cup \{0\},$ we have
$$\frac{(-1)^n(-\eta)_n}{n!}=\binom{\eta}{n}=\frac{\Gamma(\eta+1)}{n! \,\Gamma(\eta-n+1)}.$$
In particular, when $\eta=m ~(m \in \mathbb{N}_0, \, m \geq n),$ we have
$$ (-m)_n =\frac{(-1)^n \, m!}{(m-n)!}.$$

We can use this relation in Theorem \ref{jkps_T8} to obtain a family of harmonic univalent polynomials that are also
uniformly starlike in $\mathbb{D}$.

\bcor
Let $m$ be a positive real number, $c$ be a positive real numbers, $\alpha \in \mathbb{D}$ and let
$$F_1(z)=z+\overline{\frac{\alpha}{2} \sum_{n=0}^m \binom{m}{n}\frac{(m-n+1)_n}{(c)_n}z^{n+1}}, \;\;F_2(z)=z+\overline{\frac{\alpha}{2} \sum_{n=1}^m \binom{m}{n}\frac{(m-n+1)_n}{(c)_n}z^{n}}$$ and
$$ F_3(z)=z+\overline{\frac{\alpha}{2} \sum_{n=0}^m \binom{m}{n}\frac{(m-n+1)_n}{(c)_n}\frac{z^{n+1}}{n+1}}.   $$
\begin{enumerate}
\item[{\rm (a)}] If $|\alpha|\,\Gamma(c)\,\Gamma(c+2m-1)\,(m^2+2m+c-1) \leq (\Gamma(c+m))^2$,
then $F_1 \in \mathcal{US}^*_H.$\\
\item[{\rm (b)}] If $m^2|\alpha|\,\Gamma(c)\,\Gamma(c+2m-1) \leq (\Gamma(c+m))^2$,
then $F_2 \in \mathcal{US}^*_H.$\\
\item[{\rm (c)}] If $|\alpha|\,\Gamma(c)\,\Gamma(c+2m) \leq (\Gamma(c+m))^2$,
then $F_3 \in \mathcal{US}^*_H.$
\end{enumerate}
\ecor
\begin{proof}
The results follow if we set $a=b=-m$ in Theorem \ref{jkps_T8}.
\end{proof}

\noindent
{\bf Example 1}. If we let $m=3$ in Corollary 2(c), then we have the following: if $c$ is a positive real number satisfying
$$2 c + 3 c^2 + c^3 - 60 \alpha - 47 c \alpha - 12 c^2 \alpha -  c^3 \alpha \ge 0,
$$
where $\alpha$ is a complex number such that $0<|\alpha| < 1$, then the harmonic function
$$f(z)=z+\frac{\alpha}{2} \left(\overline{z}+\frac{9}{2c}\overline{z^2}+\frac{6}{c(1+c)}\overline{z^3}+\frac{3}{2c(1+c)(2+c)}\overline{z^4} \right)$$
is uniformly starlike in $\mathbb{D}.$ For instance, choosing $\alpha=e^{i \theta}/20$, we obtain that the function
$$f_{c,\,\theta}(z) = z+e^{i \theta}\left(\frac{1}{40}\overline{z}+\frac{9}{80c}\overline{z^2}+
\frac{3}{20c(1+c)}\overline{z^3}+\frac{3}{80c(1+c)(2+c)}\overline{z^4}\right)$$
is univalent and uniformly starlike in $\mathbb{D}$, whenever $c \ge 1$ and $\theta \in \IR$.
%
%
%
\bthm\label{jkps_T9}
Let either $\,a,b \in (-1, \infty)$ with $ab>0,$ or $a,b \in \mathbb{C}\setminus \{0\}$ with $b=\overline{a}$ and $\alpha \in \mathbb{D}$. Assume that $c$ is a positive real number.
\begin{enumerate}
\item[{\rm (a)}] If $c>{\rm Re\,}(a+b)+2$ and
\begin{eqnarray}\label{hyper11}
|\alpha|\frac{\Gamma(c) \Gamma(c-a-b)}{\Gamma(c-a)\Gamma(c-b)}\left[\frac{2(a)_2 (b)_2}{(c-a-b-2)_2}+\frac{5ab}{c-a-b-1}+1 \right] \leq 2,
\end{eqnarray}
then $f_1(z) \in \mathcal{UK}_H.$
\item[{\rm (b)}] If $c>{\rm Re\,}(a+b)+2$ and
\begin{eqnarray}\label{hyper12}
ab|\alpha|\frac{\Gamma(c) \Gamma(c-a-b-1)}{\Gamma(c-a)\Gamma(c-b)} \left(\frac{a+b+c+2ab}{c-a-b-2}\right)\leq 2,
\end{eqnarray}
then $f_2(z) \in \mathcal{UK}_H.$
\item[{\rm (c)}] If $c>{\rm Re\,}(a+b)+1$ and
\begin{eqnarray}\label{hyper13}
|\alpha|\frac{\Gamma(c) \Gamma(c-a-b)}{\Gamma(c-a)\Gamma(c-b)} \left(\frac{2ab+c-a-b-1}{c-a-b-1}\right)\leq 2,
\end{eqnarray}
then $f_3(z) \in \mathcal{UK}_H.$
\end{enumerate}
Here $f_1$, $f_2$ and $f_3$ are defined as in Theorem \ref{jkps_T8}.
\ethm
\begin{proof}
Following the procedure of the proof of Theorem \ref{jkps_T8}, it is enough to prove that
the coefficients of the functions $g_1(z)$, $g_2(z)$ and $g_3(z)$ defined in the proof of
Theorem \ref{jkps_T8} satisfy the sufficient condition given in Theorem \Ref{jkps_T5} for the
functions $f_k$'s ($k=1,2,3$) to be in $\mathcal{UK}_{H}$. Thus, a calculation shows that
$$K_1:= \sum_{n=1}^{\infty} n(2n-1) |b_n| \leq 1
$$
holds under the condition \eqref{hyper11}, \eqref{hyper12} and \eqref{hyper13}, respectively. Here $b_n$'s are as in the proof of Theorem \ref{jkps_T8}.
 Since the proof is similar to the proof of Theorem \ref{jkps_T8}, we omit the details here.
\end{proof}

\subsection*{Acknowledgements}
The research was supported by the project RUS/RFBR/P-163 under Department of Science \& Technology (India) and the Russian Foundation for Basic Research (project 14-01-92692). The work of the second author is supported by NBHM (DAE), India. The third author is also supported by Russian Foundation for Basic Research (Project 14-01-00510) and the Strategic Development Program of Petrozavodsk State University.

\end{document}